\input ./style/arxiv-general.cfg
\documentclass[MSNbibl,number,citesort,dvips]{arxbj}
\makeatletter
   \@ifpackageloaded{graphicx}{}{\usepackage{graphicx}}
\makeatother
\usepackage{upgreek}

\aid{0}
\volume{22}
\issue{1}
\pubyear{2016}
\firstpage{275}
\lastpage{301}
\doi{10.3150/14-BEJ658} 
\docsubty{FLA}

\makeatletter
\newtheorem{theorem}{Theorem}
\newtheorem{corollary}[theorem]{Corollary}
\newtheorem{lemma}[theorem]{Lemma}
\newtheorem{proposition}[theorem]{Proposition}
\newremark{remark}[theorem]{Remark}
\makeatother

\begin{document}
\begin{frontmatter}

\title{Quenched limit theorems for Fourier transforms and periodogram}
\runtitle{Quenched limit theorems for Fourier transforms}

\begin{aug}
\author[A]{\inits{D.}\fnms{David}~\snm{Barrera}\thanksref{e1}\ead[label=e1,mark]{barrerjd@mail.uc.edu}}%
\and
\author[A]{\inits{M.}\fnms{Magda}~\snm{Peligrad}\corref{}\thanksref{e2}\ead[label=e2,mark]{peligrm@ucmail.uc.edu}}
\address[A]{Department of Mathematical Sciences, University of Cincinnati, P.O. Box 210025,
Cincinnati, OH 45221-0025, USA. \printead{e1,e2}}
\end{aug}

\received{\smonth{1} \syear{2014}}
\revised{\smonth{4} \syear{2014}}

\begin{abstract}
In this paper, we study the quenched central limit theorem for the discrete
Fourier transform. We show that the Fourier transform of  a stationary
ergodic process, suitable centered and normalized, satisfies the quenched CLT
conditioned by the past sigma algebra. For functions of Markov chains with
stationary transitions, this means that the CLT holds with respect to the law
of the chain started at a point for almost all starting points. It is
necessary to emphasize that no assumption of irreducibility with respect to a
measure or other regularity conditions are imposed for this result. We also
discuss necessary and sufficient conditions for the validity of quenched CLT
without centering. The results are highly relevant for the study of the
periodogram of a Markov process with stationary transitions which does not
start from equilibrium. The proofs are based of a nice blend of harmonic
analysis, theory of stationary processes, martingale approximation and ergodic theory.
\end{abstract}

\end{frontmatter}

\section{Introduction}\label{s1}

The finite Fourier transform, defined as
\begin{equation}
\label{Four} S_{n}(t)=\sum_{k=1}^{n}
\mathrm{e}^{\mathrm{i}kt}X_{k},
\end{equation}
where $\mathrm{i}=\sqrt{-1}$ is the imaginary unit, plays an essential role for the
study of stationary time series $(X_{j})_{j\in\mathbb{Z}}$ of centered random
variables with finite second moment, defined on a probability space
$(\Omega,\mathcal{K},{\mathbb{P}})$.

The periodogram, introduced as a tool by Schuster \cite{Shu} in 1898, is
essential in the estimation of the spectral density of the stationary
processes. It is defined by
\begin{equation}
\label{Per} I_{n}(t)=\frac{1}{2\uppi n} \bigl|S_{n}(t)
\bigr|^{2},\qquad t\in [0,2\uppi].
\end{equation}

Wiener and Wintner \cite{Win} showed that for any stationary sequence
$(X_{j})_{j\in\mathbb{Z}}$ in ${\mathbb{L}}^{1}$ (namely $\mathbb{E}|X_{0}|<\infty$) there is a set $\Omega^{\prime}\subset  \Omega$ of
probability one such that for any $t\in [0,2\uppi]$ and any $\omega
\in\Omega^{\prime}$,  $S_{n}(t)/n$ converges. The speed of this convergence
(see Peligrad and Wu \cite{PeWu} and the references therein) is usually given
by a central limit theorem for the real and imaginary parts of $S_{n}(t)/\sqrt{n}$ under various dependence restrictions. Peligrad and Wu
\cite{PeWu} showed that, under a very mild regularity condition and finite
second moment, $[\operatorname{Re}(S_{n}(t))/\sqrt{n}, \operatorname{Im}(S_{n}(t))/\sqrt{n}]$ are asymptotically independent normal random variables
with mean $0$ and variance $\uppi f(t)$, for almost all $t$. Here $f$ is the
spectral density of $(X_{j})_{j\in\mathbb{Z}}$. This result implies that for
almost all $t$, the periodogram $I_{n}(t)$ converges in distribution to
$f(t)\chi^{2}$ where $\chi^{2}$ has a chi-square distribution with
$2$ degrees of freedom. Sufficient conditions for the validity of the law of
iterated logarithm were recently pointed out in Cuny \textit{et al.} \cite{CMP}.

An interesting problem with practical applications, is to study the validity
of the central limit theorem for Fourier transform and the periodogram for a
Markov chain with stationary transitions which is not started from equilibrium
but rather started from a point for almost all starting points. This is often
the case for simulated data and also for evolutions in random media or
particle systems. The problem is difficult, since the Markov chain started at
a point is no longer stationary. This type of central limit theorem, started
at a point, is known under the name of quenched central limit theorem (CLT)
and it is a consequence of a more general result, the almost sure conditional
CLT for stationary processes. This means that on a set of measure one the
central limit theorem holds when, in the definition of weak convergence, we
replace the usual expectation by the conditional expectation with respect to
the past $\sigma$-algebra. The almost sure conditional CLT implies CLT. Some
examples of stationary processes satisfying the CLT but not the almost sure
conditional CLT can be found in Voln\'{y} and Woodroofe \cite{VW}.

The problem of the quenched CLT for stationary Markov chains or for stationary
processes is a subject of intense research. We mention the papers \cite{DL,CV,VW,CP,CM,RS1,RS2,DMP},
among many others. Several of these results were surveyed in \cite{P}.

As far as we know, this type of convergence was not yet investigated for the
Fourier transforms or the periodograms. In this paper, we show that the
quenched CLT holds for almost all frequencies of the properly centered and
normalized discrete Fourier transform of any stationary and ergodic sequence.
We also provide necessary and sufficient conditions for the validity of
quenched CLT without centering and  specify a sufficient condition for the
validity of quenched CLT for fixed frequencies.

All these results shed additional light on the speed of convergence of the
periodogram in approximating the spectral density $f(t)$ of a stationary
process. The techniques are a nice blend of martingale approximation, rooted
in Gordin \cite{G1} and Rootz\'{e}n \cite{Rootzen} and developed by Gordin and
Lif{\v{s}ic} \cite{GL} and Woodroofe \cite{Wood}, and tools from ergodic theory
and harmonic analysis.

To allow for flexibility in applications, we introduce a stationary sequence
and a filtration in two different ways. First by using a measure preserving
transformation, and then, in Section~\ref{s3}, as a function of a Markov chain. We
formulate the main results for measure preserving transformations in terms of
almost sure conditional CLT. However, in Section~\ref{s3} we show that, only by a
change of language, the results can be formulated for stationary and ergodic
Markov chains, where the terminology of a process started at a point becomes natural.

A variety of applications to functions of linear processes, functions of
Markov chains, iterated random functions, mixing sequences, are also pointed
out. It is remarkable that for the case of a stationary ergodic reversible
Markov chain the quenched CLT without centering holds without any other
additional assumptions.

Our paper is organized as follows. Section~\ref{s2} contains the presentation of the
results. Several applications are given in Section~\ref{s3}. Section~\ref{s4} is devoted to
the proofs. Section~\ref{s5} contains several auxiliary results needed for the main proofs.

\section{Definitions, background and results}\label{s2}

A strictly stationary sequence can be introduced in many equivalent ways. It
can be viewed, for instance, as a function of a stationary Markov chain with
general state space. This definition will be given in Section~\ref{s3}. For more
flexibility in the selection of filtration, in this section, we shall
introduce a stationary sequence and a filtration by using a measure preserving transformation.

Let $(\Omega,\mathcal{K},{\mathbb{P}})$ be a probability space where, without
restricting the generality, we shall assume that $\mathcal{K}$ is countably
generated, and let $T\dvtx \Omega\mapsto\Omega$ be a bijective bi-measurable
transformation preserving the probability ${\mathbb{P}}$. An element $A$ is
said to be invariant if $T(A)=A$. We denote by $\mathcal{I}$ the $\sigma$-algebra of all invariant sets.
The transformation $T$ is ergodic with
respect to $\mathbb{P}$ if each element of $\mathcal{I}$ has probability $0$
or $1$. Let $\mathcal{F}_{0}$ be a $\sigma$-algebra of $\mathcal{K}$
satisfying $\mathcal{F}_{0}\subseteq T^{-1}(\mathcal{F}_{0})$. Define the
nondecreasing filtration $(\mathcal{F}_{i})_{i\in{\mathbb{Z}}}$ by
$\mathcal{F}_{i}=T^{-i}(\mathcal{F}_{0})$ and let ${\mathcal{F}}_{-\infty
}=\bigcap_{k\in{\mathbb{Z}}}{\mathcal{F}}_{k}$. Let $X_{0}$ be a
$\mathcal{F}_{0}$-measurable, square integrable and centered random variable.
Define the sequence $\mathbf{X}=(X_{i})_{i\in\mathbb{Z}}$ by
\begin{equation}
\label{defX2} X_{i}=X_{0}\circ T^{i}.
\end{equation}
For $p\geq 1$,  we denote by  $\Vert\cdot\Vert_{p}$ the norm in $\mathbb{L}^{p}(\Omega,\mathcal{F},\mathbb{P})$
and for an integrable random variable $Y$
we denote by $\mathbb{E}_{0}(Y)=\mathbb{E}(Y|\mathcal{F}_{0})$.

Since $\mathcal{K}$ is countably generated, there is a regular conditional
probability measure $\mathbb{P}^{\omega}(\cdot)$ with respect to
$\mathcal{F}_{0}$, such that for all $\omega\in\Omega$, $\mathbb{P}^{\omega}(\cdot)$
is a measure on $\mathcal{K}$ and for each $A\in\mathcal{K}$ we have
$\mathbb{P}^{\omega}(A)=\mathbb{P}(A|\mathcal{F}_{0})(\omega)$, $\mathbb{P}$
a.s. For integrable $X$, the corresponding conditional expectation is
denoted by $\mathbb{E}^{\omega}(X)$ and it is a regular version of
$\mathbb{E}(X|\mathcal{F}_{0})(\omega)$.

Relevant to our results is the notion of spectral distribution function
induced by the covariances. By Herglotz's theorem (see, e.g., Brockwell and
Davis \cite{BD}), there exists a nondecreasing function $G$ (the spectral
distribution function) on $[0,2\uppi]$ such that, for all $j\in\mathbb{Z}$,
\[
\operatorname{cov}(X_{0},X_{j})=\int_{0}^{2\uppi}
\mathrm{e}^{\mathrm{i}j\theta}\,\mathrm{d}G(\theta),\qquad j\in\mathbb{Z}.
\]
If $G$ is absolutely continuous with respect to the normalized Lebesgue
measure $\lambda$ on $[0,2\uppi]$, then the Radon--Nikodym derivative $f$ of $G$
with respect to the Lebesgue measure is called the spectral density and we
have
\[
\operatorname{cov}(X_{0},X_{j})=\int_{0}^{2\uppi}
\mathrm{e}^{\mathrm{i}j\theta}f(\theta )\,\mathrm{d}\theta, \qquad j\in\mathbb{Z}.
\]
We shall introduce the notations
\[
\mathbf{V}_{n}(t)=\frac{1}{\sqrt{n}} \bigl[\operatorname{Re}
\bigl(S_{n}(t) \bigr), \operatorname{Im} \bigl(S_{n}(t) \bigr)
\bigr].
\]
We also denote
\[
\mathbf{W}_{n}(t)=\frac{1}{\sqrt{n}} \bigl[\operatorname{Re}
\bigl(S_{n}(t)-\mathbb{E}_{0} S_{n}(t) \bigr),
\operatorname{Im} \bigl(S_{n}(t)-\mathbb{E}_{0}S_{n}(t)
\bigr) \bigr].
\]
The central limit theorem for $\mathbf{V}_{n}(t)$ has a long history. We
mention, among many others, Rosenblatt (Theorem~5.3, page~131, \cite{ro85}) who
considered mixing processes; Brockwell and Davis (Theorem~10.3.2, page~347,
\cite{BD}), Walker \cite{Walker} and Terrin--Hurvich \cite{Ter} discussed
linear processes; Wu \cite{Wu05} treated mixingales.

Peligrad and Wu \cite{PeWu} established the following result, where, besides a
mild regularity assumption (\ref{regular}), no other restriction of dependence
is imposed to the stochastic process. Below, by~$\Rightarrow$ we denote
convergence in distribution.

\renewcommand{\thetheorem}{\Alph{theorem}}
\begin{theorem}[(Peligrad and Wu)]\label{thm1}
Let  $(X_{k})_{k\in\mathbb{Z}}$ be a stationary ergodic process, centered, with finite second
moments, such that the following regularity assumption is satisfied,
\begin{equation}
\label{regular} \mathbb{E}(X_{0}|\mathcal{F}_{-\infty})=0\qquad
\mathbb{P}\mbox{-a.s.}
\end{equation}
Then, for almost all $t\in(0,2\uppi)$, the following
convergence holds:
\begin{equation}
\label{defg} \lim_{n\rightarrow\infty}{\frac{\mathbb{E}{|S_{n}(t)|^{2}}}{n}}=2\uppi f(t),
\end{equation}
where  $f(t)$ is the spectral density of $(X_{k})_{k\in\mathbb{Z}}$.
Furthermore
\begin{equation}
\label{ACLT} \frac{1}{\sqrt{n}}\mathbf{V}_{n}(t)\Rightarrow
\mathbf{N}(t)\qquad \mbox{under }\mathbb{P},
\end{equation}
where  $\mathbf{N}(t)=[N_{1}(t), N_{2}(t)]$,  with $N_{1}(t)$ and  $N_{2}(t)$ independent identically distributed normal random
variables mean  $0$ and variance  $\uppi f(t)$.
\end{theorem}

The proof of Theorem~\ref{thm1}  is based on the celebrated Carleson's \cite{Ca} theorem
on almost sure convergence of Fourier transforms. A different proof, without
using Carleson's result, was recently given in Cohen and Conze \cite{CC}. This
suggests that the power of Carleson's theorem might lead to a stronger type of
limiting distribution, in the almost sure sense.

The goal of this paper is to study a more general form of (\ref{ACLT}) known
under the name of quenched CLT.

By the quenched CLT we shall understand the following almost sure conditional
limit theorem:

For almost all $t\in [0,2\uppi]$, there is $\Omega^{\prime}$ with
$\mathbb{P}(\Omega^{\prime})=1$ such that for all $\omega\in\Omega^{\prime}$
we have
\begin{equation}
\label{QCLT1}
\mathbb{E}^{\omega}\bigl[g \bigl(\mathbf{V}_{n}(t) \bigr)\bigr]
\rightarrow \mathbb{E} \bigl[g \bigl(\mathbf{N}(t) \bigr) \bigr]\qquad\mbox{as }n
\rightarrow\infty,
\end{equation}
for any function $g$ which is continuous and bounded. We shall say in this
case that the quenched CLT holds for almost all frequencies. In other
notation, for almost all $t\in [0,2\uppi]$, there is $\Omega^{\prime}$ with
$\mathbb{P}(\Omega^{\prime})=1$ such that for all $\omega\in\Omega^{\prime}$
we have
\[
\mathbf{V}_{n}(t)\Rightarrow\mathbf{N}(t)\qquad \mbox{as } n
\rightarrow\infty\mbox{ under }\mathbb{P}^{\omega}.
\]
Clearly (\ref{QCLT1}) implies (\ref{ACLT}) by integration with respect to
$\mathbb{P}$.

Our first result gives a quenched CLT under a certain centralization. Note
that the next theorem applies to any stationary and ergodic sequence.

\renewcommand{\thetheorem}{\arabic{theorem}}
\setcounter{theorem}{0}
\begin{theorem}
\label{c-quenched}
Let $(X_{k})_{k\in\mathbb{Z}}$ be a stationary ergodic
process, $X_{k}$ defined by (\ref{defX2}) and let ${S_{n}(t)}$ be defined by
(\ref{Four}). Then, for almost all $t\in [0,2\uppi]$
\[
\lim_{n\rightarrow\infty}{\frac{\mathbb{E}_{0}{|S_{n}(t)-\mathbb{E}}_{0}%
{S_{n}(t)|}^{2}}{n}}=\sigma_{t}^{2}\qquad
\mathbb{P}\mbox{-a.s.}
\]
and the quenched CLT holds for $\mathbf{W}_{n}(t)$, where $N_{1}(t)$ and
$N_{2}(t)$ are independent identically distributed normal random
variables with mean $0$ and variance $\sigma_{t}^{2}/2$.
\end{theorem}

Our second theorem provides a characterization of quenched convergence without centering.

\begin{theorem}
\label{char}
Let $(X_{k})_{k\in\mathbb{Z}}$ be as in Theorem~\ref{c-quenched}.
Then the following statements are equivalent:
\begin{enumerate}[(a)]
\item[(a)] For almost all $t\in [0,2\uppi]$ the quenched CLT in (\ref{QCLT1}) holds,
where $N_{1}(t)$ and $N_{2}(t)$ are as in Theorem~\ref{c-quenched}.
\item[(b)] For almost all $t\in [0,2\uppi]$ we have $\frac{1}{\sqrt{n}}\mathbb{E}_{0}(S_{n}(t))\rightarrow0$,
$\mathbb{P}$-a.s.
\end{enumerate}
\end{theorem}

\subsection*{Discussion}
An interesting problem is to specify $\sigma_{t}^{2}$.
Note first that if $\mathbb{E}_{0}(S_{n}(t))/\sqrt{n}$ converges in
$\mathbb{L}^{2}$ to $0$, our proofs show that in this case $\sigma_{t}^{2}$
can be identified as
\begin{equation}
\label{defsi} \lim_{n\rightarrow\infty}\frac{\mathbb{E}{|S_{n}(t)|^{2}}}{n}=
\sigma_{t}^{2}.
\end{equation}
Note also that in both Theorems~\ref{c-quenched} and~\ref{char} we do not
require the sequence to be regular, that is, it may happen that $E(X_{0}|\mathcal{F}_{-n})$ does not converge to $0$ in $\mathbb{L}^{2}$. The spectral
density might not exist. If we assume condition (\ref{regular}) then, as shown
in Peligrad and Wu \cite{PeWu}, the spectral density $f(t)$ of $(X_{k})_{k\in\mathbb{Z}}$ exists and $\sigma_{t}^{2}=2\uppi f(t)$. Furthermore, we have
\[
\frac{1}{\sqrt{n}}\mathbb{E}_{0} \bigl(S_{n}(t) \bigr)
\rightarrow 0 \qquad \mbox{in } \mathbb{L}^{2},
\]
and then $\sigma_{t}^{2}$ can also be identified as
\begin{equation}
\label{defsimat} \lim_{n\rightarrow\infty}{\frac{\mathbb{E}{|S_{n}(t)|^{2}}}{n}}=
\sigma_{t} ^{2}=2 \uppi f(t).
\end{equation}

By the mapping theorem (see Theorem~29.2 in \cite{BilPrM}), all our results
imply corresponding results for the periodogram. As a consequence of Theorem~\ref{char} and the discussion above we obtain the following corollary:

\begin{corollary}
Assume that the sequence $(X_{k})_{k\in\mathbb{Z}}$ is as in Theorem~\ref{c-quenched}
and in addition satisfies (\ref{regular}) and item
\textup{(b)} of Theorem~\ref{char}. Then, for almost all $t\in [0,2\uppi]$ the
periodogram $I_{n}=(2\uppi n)^{-1}|S_{n}(t)|^{2}$ satisfies a quenched limit
theorem with the limit $f(t)\chi^{2}(2)$, where $\chi^{2}(2)$ is a chi-square
random variable with $2$ degrees of freedom, and $f(t)$ the spectral density.
\end{corollary}

The next corollary provides sufficient conditions for the validity of item (b)
of Theorem~\ref{char}.

\begin{corollary}
\label{corquenched}
Assume that $(X_{k})_{k\in\mathbb{Z}}$ is as in Theorem~\ref{c-quenched} and in addition that
\begin{equation}
\label{condcorquenched} \sum_{k\geq1}\frac{|\mathbb{E}_{0}(X_{k+1}-X_{k})|^{2}}{k}<\infty
\qquad \mathbb{P}\mbox{-a.s.}
\end{equation}
Then (\ref{QCLT1}) holds, where $N_{1}(t)$ and $N_{2}(t)$ are as in Theorem~\ref{c-quenched}.
\end{corollary}

Clearly (\ref{condcorquenched}) is satisfied if
\begin{equation}
\label{sufcond0} \sum_{k\geq1}\frac{|\mathbb{E}_{0}(X_{k})|^{2}}{k}<\infty
\qquad \mathbb{P}\mbox{-a.s.},
\end{equation}
which is further implied by
\begin{equation}
\label{sufcond}%
\sum_{k\geq1}\frac{\Vert\mathbb{E}_{0}(X_{k})\Vert_{2}^{2}}{k}<
\infty.
\end{equation}
Moreover, since $\Vert\mathbb{E}_{0}(X_{k})\Vert_{2}$ is decreasing, condition
(\ref{sufcond}) implies condition (\ref{regular}). These remarks justify the
following corollary:

\begin{corollary}
\label{corQuencedf}
Condition (\ref{sufcond}) is sufficient for the quenched
CLT in (\ref{QCLT1}) with $N_{1}(t)$ and $N_{2}(t)$ i.i.d. normal
random variables with mean  $0$ and variance $\uppi f(t)$, $f(t)$ being
the spectral density of the process.
\end{corollary}

The above results hold for almost all frequencies. Actually it is possible
that on a set of measure $0$ the behavior be quite different. For the case
when $t=0$, there are a variety of examples where the partial sums of a
stationary sequence do not satisfy a nondegenerate CLT. One important example
of this kind is provided by filters of Gaussian processes with long range
dependence, when the covariances are not summable. For example, Rosenblatt
\cite{rosen61} proved that for a stationary Gaussian sequence $(X_{k})_{k\in\mathbb{Z}}$
of standard normal random variables with $\operatorname{cov}(X_{0},X_{k})=(1+k^{2})^{-\alpha/2}, 0<\alpha<1/2$,
the sequence $n^{-1+\alpha
}\sum_{k=1}^{n}(X_{k}^{2}-1)$ has a nonnormal limiting distribution as
$n\rightarrow\infty$. Another interesting example, also for $t=0$, is provided
by Herrndorf \cite{Hern} who constructed a stationary sequence of centered
uncorrelated random variables with finite second moment, which is strongly
mixing with arbitrary mixing rate and the partial sums do not satisfy a
nondegenerate CLT under any normalization converging to infinite. This
example satisfies condition (\ref{regular}). Furthermore, Bradley
\cite{Bradley} (see Theorem~34.14, Vol.~3) constructed a stationary sequence
of centered random variables with finite second moment, satisfying our
condition (\ref{sufcond}) and such that its partial sums normalized by its
standard deviation is attracted to a non-Gaussian nondegenerate distribution.
Rosenblatt \cite{Ro81} studied the Fourier transform of nonlinear functions of
Gaussian processes and established for certain frequencies, on a set of
measure $0$, non-Gaussian attraction for the Fourier transform properly normalized.

In the spirit of Maxwell and Woodroofe \cite{MW} and Cuny and Merlev\`{e}de
\cite{CM}, we give below a result allowing us to identify frequencies for
which the quenched CLT holds.

\begin{theorem}
\label{MW}
Let $t\in(0,2\uppi)$ be  such that $\mathrm{e}^{-2\mathrm{i}t}$ is not an
eigenvalue of $T$. Assume that the sequence $(X_{k})_{k\in\mathbb{Z}}$ is as
in Theorem~\ref{c-quenched}  and in addition that we have
\begin{equation}
\label{condMW} \sum_{k\geq1}\frac{1}{k^{3/2}} \bigl\Vert
\mathbb{E}_{0} \bigl(S_{k}(t) \bigr) \bigr\Vert_{2}<
\infty.
\end{equation}
Then (\ref{QCLT1}) holds with $N_{1}(t)$ and $N_{2}(t)$ independent
identically distributed normal random variables mean  $0$ and variance
$\sigma_{t}^{2}/2$ where $\sigma_{t}^{2}$ is identified by (\ref{defsi}).
\end{theorem}

\section{Applications}\label{s3}

\subsection{Functions of Markov chains}\label{s31}

Let $(\xi_{n})_{n\in\mathbb{Z}}$ be a stationary and ergodic Markov chain
defined on a probability space $(\Omega,\mathcal{F},\mathbb{P})$ with values
in a measurable space $(S,\mathcal{A})$. The marginal distribution is denoted
by $\pi(A)=\mathbb{P}(\xi_{0}\in A)$ and we assume that there is a regular
conditional distribution for $\xi_{1}$ given $\xi_{0}$ denoted by
$Q(x,A)=\mathbb{P}(\xi_{1}\in A| \xi_{0}=x)$. In addition $Q$ denotes the
Markov operator acting via $(Qh)(x)=\int_{S}h(s)Q(x,\mathrm{d}s)$. Next, let
$\mathbb{L}_{0}^{2}(\pi)$ be the set of measurable functions on $S$ such that
$\int h^{2}\,\mathrm{d}\pi<\infty$ and $\int h\,\mathrm{d}\pi=0$. For a function
$h\in \mathbb{L}_{0}^{2}(\pi)$ let
\begin{equation}
\label{defX} {X_{i}=h(\xi_{i})}.
\end{equation}
Denote by $\mathcal{F}_{k}$ the $\sigma$-field generated by $\xi_{i}$ with
$i\leq k$. For any integrable random variable $X$ we denote $\mathbb{E}_{k}(X)=\mathbb{E}(X|\mathcal{F}_{k})$ and
$\mathbb{P}_{k}(A)=\mathbb{P}(A|\mathcal{F}_{k})$.
In our notation $\mathbb{E}_{0}(X_{1})=(Qh)(\xi
_{0})=\mathbb{E}(X_{1}|\xi_{0})$.

To guarantee that the regular transitions exist, we shall assume that
$\mathcal{A}$ is countably generated.

The Markov chain is usually constructed in a canonical way on $\Omega
=S^{\infty}$ endowed with sigma algebra $\mathcal{A}^{\infty}$, and
$\mathcal{\xi}_{n}$ is the $n$th projection on $S$. The shift
$T\dvtx \Omega\rightarrow\Omega$ is defined by $\mathcal{\xi}_{n}(T\omega
)=\mathcal{\xi}_{n+1}(\omega)$ for every $n\geq0$.

For any probability measure $\upsilon$ on $\mathcal{A}$ the law of $(\xi
_{n})_{n\in\mathbb{Z}}$ with transition operator $Q$ and initial distribution
$\upsilon$ is the probability measure $\mathbb{P}^{\upsilon}$ on $(S^{\infty
},\mathcal{A}^{\infty})$ such that
\[
\mathbb{P}^{\upsilon}(\xi_{n+1}\in A| \xi_{n}=x)=Q(x,A)
\quad\mbox{and}\quad \mathbb{P}^{\upsilon}(\xi_{0}\in A)=
\upsilon(A).
\]
For $\upsilon=\pi$,  we denote $\mathbb{P}=\mathbb{P}^{\pi}$.
For $\upsilon
=\delta_{x}$, the Dirac measure, denote by $\mathbb{P}^{x}$ and $\mathbb{E}^{x}$
the regular probability and conditional expectation for the process
started at $x$. Note that for each $x$ fixed $\mathbb{P}^{x}(\cdot)$ is a
measure on $\mathcal{F}^{\infty}$, the sigma algebra generated by
$\bigcup_{k}\mathcal{F}_{k}$. Furthermore $\mathbb{P}^{x}(\cdot)$ is a version of the
conditional probability on $\mathcal{F}^{\infty}$ given $\xi_{0}$ and, by
Markov property, $\mathbb{P}^{x}(\cdot)$ is also the regular measure on
$\mathcal{F}^{\infty}$ given $\mathcal{F}_{0}$.

We mention that any stationary sequence $(Y_{k})_{k\in\mathbb{Z}}$ can be
viewed as a function of a Markov process $\xi_{k}=(Y_{j};j\leq k$) with the
function $g(\xi_{k})=Y_{k}$. Therefore the theory of stationary processes can
be embedded in the theory of Markov chains.

For a Markov chain, by the quenched CLT for the Fourier transform we shall
understand the following convergence: for almost all $t\in [0,2\uppi]$
there is a set $S^{\prime}\subset S$ with $\pi(S^{\prime})=1$ such that for
$x\in S^{\prime}$
\begin{equation}
\label{QCLT} \mathbf{V}_{n}(t)\Rightarrow\mathbf{N}(t)\qquad
\mbox{under }\mathbb{P}^{x}.
\end{equation}
In other words for almost all $t\in [0,2\uppi]$, there is a set $S^{\prime
}\subset S$ with $\pi(S^{\prime})=1$ such that for $x\in S^{\prime}$
\[
\mathbb{E}^{x}\bigl[g \bigl(\mathbf{V}_{n}(t) \bigr)\bigr]\rightarrow
\mathbb{E}\bigl[g \bigl(\mathbf{N}(t) \bigr)\bigr]\qquad\mbox{as }n\rightarrow\infty,
\]
for any function $g$ continuous and bounded. When the stationary process is
viewed as a function of Markov chain, then $\xi_{0}=(Y_{j};j\leq0)$, and
therefore a fixed value of $\xi_{0}$ means a fixed past trajectory up to the
moment of time $0$.

All our results hold in the setting of Markov chains. In this case, the
transformation $T$ is the shift. The Markov property allows for the
formulation (\ref{QCLT}).

It is remarkable that for ergodic reversible Markov chains the quenched CLT
holds without centering and without any additional assumptions.

\begin{corollary}
\label{revcopy1}
Assume that $(X_{k})_{k\in\mathbb{Z}}$ is defined by
(\ref{defX}) and in addition that the Markov chain $(\xi_{k})_{k\in Z}$ is
reversible (i.e. $Q=Q^{\ast})$. Let $t\in(0,2\uppi)\setminus\{\uppi, \uppi/2, 3\uppi/2\}$.
Then, (\ref{QCLT}) holds where $N_{1}(t)$ and $N_{2}(t)$ are as in Theorem~\ref{MW}.
\end{corollary}

\begin{pf}
We shall verify the conditions of Theorem~\ref{MW}. Since the spectrum of $Q$
is contained in $[-1,1]$ and for $t\in(0,2\uppi)\setminus\{\uppi\}$ we have that
$\mathrm{e}^{\mathrm{i}t}$ is not real, the operator $I-\mathrm{e}^{\mathrm{i}t}Q$ is
invertible, and therefore there exists $g\in{\mathbb{L}}^{2}(S,\mathcal{A},\pi)$ such that $h=g-\mathrm{e}^{\mathrm{i}t}Qg$. We obtain
\[
\mathbb{E}_{0} \bigl(S_{n}(t) \bigr)=\sum
_{k=1}^{n}\mathbb{E}_{0} \bigl(
\mathrm{e}^{\mathrm{i}tk} g(\xi_{k})-\mathrm{e}^{\mathrm{i}t(k+1)}g(
\xi_{k+1}) \bigr)=\mathrm{e}^{\mathrm{i}t}\mathbb{E} _{0}
\bigl(g(\xi_{1}) \bigr)-\mathrm{e}^{\mathrm{i}t(n+1)}\mathbb{E}_{0}
\bigl(g(\xi_{n+1}) \bigr).
\]
Then clearly
\[
\bigl\Vert\mathbb{E}_{0} \bigl(S_{n}(t) \bigr)
\bigr\Vert_{2} \leq2 \Vert g\Vert_{2}<\infty,
\]
and therefore condition (\ref{condMW}) is satisfied.
Furthermore,  since $T$ is the shift operator, under our hypotheses, cannot
have eigenvalues other than ${\pm}1$ (see page 15 in Cuny \textit{et al.} \cite{CMP}).
\end{pf}

\subsection{Iterated random functions}\label{s32}

Let $(\Gamma,d)$ be a complete and separable metric space and let $\xi
_{n}=F_{\varepsilon_{n}}(\xi_{n-1})$, where $F_{\varepsilon}(\cdot
)=F(\cdot,\varepsilon)$ is the $\varepsilon$-section of a jointly measurable
function $F\dvtx \Gamma\times\Upsilon \rightarrow\Gamma$ and $\varepsilon
,\varepsilon_{n}$,  $n\in Z$ are i.i.d. random variables taking values in a
second measurable space $\Upsilon$. Define $L_{\varepsilon}=\sup_{x\neq
x^{\prime}}d(F_{\varepsilon}(x);F_{\varepsilon}(x^{\prime}))/d(x,x^{\prime
})$. Diaconis and Freedman \cite{DiFr} proved that $(\xi_{n})$ admits a unique
stationary distribution $\pi$ provided that for some $\alpha>0$ and $x_{0}\in\Gamma$,
\begin{equation}
\label{alphacond} \mathbb{E} \bigl(L_{\varepsilon}^{\alpha} \bigr)<\infty,
\qquad \mathbb{E}(\log L_{\varepsilon})<0\quad \mbox{and}\quad \mathbb{E}
\bigl(d^{\alpha} \bigl(x_{0},F_{\varepsilon
}(x) \bigr) \bigr)<
\infty.
\end{equation}
Let $h$ be a function and let $X_{k}=h(\xi_{k})$.
Assume $\mathbb{E}(X_{1})=0$
and $\mathbb{E}|X_{1}|^{2}<\infty$.
To analyze this example, we shall use the
coupling function introduced by Wu \cite{Wu05}:
\[
\Delta_{h}(t)=\sup \bigl\Vert \bigl(h(\xi)-h \bigl(\xi^{\prime}
\bigr) \bigr)I \bigl(d \bigl(\xi,\xi^{\prime} \bigr)<t \bigr)
\bigr\Vert_{2},
\]
where the supremum is taken over all $\xi,\xi^{\prime}$ independent
distributed as $\pi$. We shall establish the following:

\begin{corollary}
Assume condition (\ref{alphacond}) is satisfied and
\begin{equation}
\label{deltacond} \int_{0}^{1/2}\frac{\Delta_{h}^{2}(t)}{t|\log t|}\,
\mathrm{d}t< \infty.
\end{equation}
Then, for almost all frequencies, the quenched CLT (\ref{QCLT}) holds with
$N_{1}(t)$ and $N_{2}(t)$ i.i.d.  normal random variables with mean
$0$  and variance $\uppi f(t)$, $f(t)$ being the spectral density of
the process.
\end{corollary}

\begin{pf}
We shall verify condition (\ref{sufcond0}). By Lemma~3 in Wu
and Woodroofe \cite{WuWo}, condition (\ref{alphacond}) implies that there is
$\beta>0$, $C>0$ and $0<r<1$ such that
\begin{equation}
\label{ww} \mathbb{E} \bigl(d^{\beta} \bigl(\xi_{n},
\xi_{n}^{\prime} \bigr) \bigr)\leq Cr^{n},
\end{equation}
where $\xi_{n},\xi_{n}^{\prime}$ are i.i.d. Since $\mathbb{E}(h(\xi
_{n}^{\prime})|\xi_{0})=0$ a.s.
\begin{eqnarray*}
\bigl|\mathbb{E} \bigl(h(\xi_{n} ) |\xi_{0} \bigr)\bigr| &\leq & \bigl|
\mathbb{E} \bigl( \bigl[h (\xi_{n} )-h \bigl(\xi_{n}^{\prime
}
\bigr) \bigr]I \bigl(d \bigl(\xi_{n},\xi_{n}^{\prime}
\bigr) \leq  \delta_{n} \bigr) |\xi_{0} \bigr) \bigr|
\\
&&{}+\bigl|\mathbb{E} \bigl( \bigl[h(\xi_{n})-h \bigl(\xi_{n}^{\prime}
\bigr) \bigr]I \bigl(d \bigl(\xi_{n},\xi_{n}^{\prime
}
\bigr)>\delta_{n} \bigr)|\xi_{0} \bigr)\bigr|\\
&=& I_{n}+
\mathit{II}_{n}.
\end{eqnarray*}
To establish (\ref{sufcond0}), it is enough to prove that
\begin{equation}
\label{one} \sum_{n\geq1}\frac{\mathbb{E}(I_{n}^{2})}{n}<\infty
\end{equation}
and%
\begin{equation}
\label{two} \sum_{n\geq1}\frac{\mathit{II}_{n}^{2}}{n}<\infty
\qquad\mbox{a.s.}
\end{equation}
By Cauchy--Schwartz inequality and Markov inequality
\[
\mathbb{E} (\mathit{II}_{n} )\leq2^{1/2} \Vert
X_{0} \Vert_{2}\mathbb{P}^{1/2} \bigl(d \bigl(
\xi_{n},\xi _{n}^{\prime} \bigr)>\delta_{n}
\bigr)\leq \Vert X_{0} \Vert_{2} \bigl[2\mathbb{E}
\bigl(d^{\beta} \bigl( \xi_{n},\xi_{n}^{\prime}
\bigr) \bigr)/ \delta_{n}^{\beta} \bigr]^{1/2}.
\]
By selecting now $\delta_{n}=r^{n/2\beta}$ we obtain $\mathbb{E}(\mathit{II}_{n})\leq
r^{n/4}$.
Therefore $\mathbb{P}(\mathit{II}_{n}>r^{n/8})\leq r^{n/4}$,
and (\ref{two})
follows by the Borel--Cantelli lemma.

Next, note that $I_{n}^{2}\leq\Delta_{h}^{2}(\delta_{n})$ and for the
selection of $\delta_{n}=r^{n/2\beta}$, the convergence of the series in
(\ref{one}) holds under the integral condition (\ref{deltacond}).

Furthermore, the above computations also show that $\mathbb{E}|\mathbb{E}(
h(\xi_{n})|\xi_{0})|\rightarrow0$ as $n\rightarrow\infty$ which proves
(\ref{regular}).
\end{pf}

\subsection{Linear processes}\label{s33}

Next, we give an application to linear processes.

\begin{corollary}
\label{linear}
Let $(\xi_{k})_{k\in\mathbb{Z}}$ be a sequence of stationary and
ergodic square integrable martingale differences. Define
\begin{equation}
\label{deflinear}
X_{k}=\sum_{j=0}^{\infty}a_{j}
\xi_{k-j},\qquad\mbox{where }\sum_{j=0}^{\infty}
a_{j}^{2}< \infty.
\end{equation}
Then, under the condition
\[
\sum_{j\geq 3}(a_{j}-a_{j+1})^{2}
\log j<\infty,
\]
the conclusion of Corollary~\ref{corquenched} holds.
\end{corollary}

\begin{pf}
We shall verify the conditions of Corollary~\ref{corquenched}.

Clearly for $k\geq 1$,  by the orthogonality of the martingale differences
\begin{eqnarray*}
\bigl\Vert\mathbb{E}_{0}(X_{k+1}-X_{k})
\bigr\Vert_{2}^{2} &=& \biggl\Vert\mathbb{E}_{0} \biggl(\sum
_{j\geq
-1}a_{j+1}\xi_{k-j}-\sum
_{j\geq0}a_{j}\xi_{k-j} \biggr)
\biggr\Vert_{2}^{2}
\\
&=& \biggl\Vert\sum_{j\geq k}a_{j+1}\xi_{k-j}-
\sum_{j\geq k}a_{j}\xi_{k-j}
\biggr\Vert_{2}^{2}
=\sum_{j\geq k}
( a_{j+1}-a_{j} ) ^{2}\Vert\xi_{0}
\Vert_{2}^{2}.
\end{eqnarray*}
Now
\[
\sum_{k\geq1}\frac{1}{k}\sum
_{j\geq k} ( a_{j+1}-a_{j} )
^{2}%
\leq\sum_{j\geq1} (
a_{j+1}-a_{j} ) ^{2}\log j,
\]
and the conclusion follows by Corollary~\ref{corquenched}.
\end{pf}

\begin{remark}
In the case when the sequence $a_{j}$ is positive and decreasing, then the
natural condition $\sum_{j=0}^{\infty}a_{j}^{2}<\infty$ is necessary and
sufficient for the conclusion of Corollary~\ref{linear}.
\end{remark}

\subsection{Functions of linear processes}\label{s34}

In this section, we shall focus on functions of real-valued linear processes.
Let $(a_{i})_{i\in\mathbb{Z}}$ be a~sequence of square summable real numbers
and $(\xi_{i})_{i\in\mathbb{Z}}$ is a sequence of i.i.d. random variables in
${\mathbb{L}}^{2}$ with mean $0$ and variance $\sigma^{2}$. Define $X_{k}$ by
(\ref{deflinear}) and let $h$ be a real valued function and define
\begin{equation}\label{def2suite}
Y_{k}=h(X_{k})-{\mathbb{E}}h(X_{k}).
\end{equation}
As in \cite{CMP} we shall give sufficient conditions for the validity of (\ref{QCLT1}) in terms of the modulus
of continuity of the function~$h$ on the interval $[-M,M]$, defined by
\begin{equation}\label{modulusbound}
w_{h}(u,M)=\sup \bigl\{\bigl|h(x)-h(y)\bigr|,|x-y|\leq u,|x|\leq M,|y|\leq M
\bigr\}.
\end{equation}

\begin{corollary}
\label{corlin}
Assume that $h$ is $\gamma$-H\"{o}lder on any compact set, with
$w_{h}(u,M)\leq Cu^{\gamma}M^{\beta}$, for some $C>0$, $\gamma\in(0,1]$ and
$\beta\geq0$. Assume that $\mathbb{E}(h^{2}(X_{k}))<\infty$ and
\begin{equation}\label{condcorflin}
\sum_{k\geq3}a_{k}^{2}\log k<
\infty\quad\mbox{and}\quad \mathbb{E} |\xi_{0}|^{2\vee2\gamma\vee2\beta}<
\infty.
\end{equation}
Then (\ref{QCLT1}) holds with $N_{1}(t)$ and $N_{2}(t)$ i.i.d. normal
random variables, mean $0$ and variance $\uppi f(t)$, $f(t)$ being the
spectral density of the process.
\end{corollary}

\begin{pf}
We shall apply Corollary~\ref{corQuencedf}. Define
$\mathcal{F}_{k}=\sigma(\xi_{l},l\leq k)$. Since $\mathcal{F}_{-\infty}$
is trivial, (\ref{regular}) holds. We write
\[
Y_{0}=\sum_{l\geq0}P_{-l}(Y_{0}),
\]
where $P_{-l}$ denotes the projector operator%
\begin{equation}\label{projop}
P_{-l}(\cdot)={\mathbb{E}}_{-l}(\cdot)-{
\mathbb{E}}_{-l-1}(\cdot).
\end{equation}
By the orthogonality of the projections,
\[
\bigl\Vert{\mathbb{E}}_{0}(Y_{k})\bigr\Vert_{2}^{2}=
\sum_{l\geq0}\bigl\Vert P_{-l}
(Y_{k})
\bigr\Vert_{2}^{2}=\sum_{j\geq k}\bigl\Vert
P_{0}(Y_{j})\bigr\Vert_{2}^{2}<\infty.
\]
Therefore, condition (\ref{sufcond}) follows from
\begin{equation}\label{condproj}
\sum_{j\geq2}\bigl\Vert P_{0}(Y_{j})
\bigr\Vert_{2}^{2}\log j<\infty.
\end{equation}
So it remains to verify (\ref{condproj}). We estimate $\Vert P_{0}(Y_{j}%
)\Vert_{2}^{2}$ as in \cite{CMP}. We give here the argument for completeness.
Let $\xi^{\prime}$ be an independent copy of $\xi$, and denote by
${\mathbb{E}}_{\xi}(\cdot)$ the conditional expectation with respect to $\xi$.
Clearly
\[
P_{0}(Y_{k})=\mathbb{E}_{\xi} \Biggl[h \Biggl(
\sum_{j=0}^{k-1}a_{j}
\xi_{k-j}^{\prime
}+a_{k}\xi_{0}+\sum
_{j>k}a_{j}\xi_{k-j} \Biggr)-h \Biggl(\sum
_{j=0}^{k-1}a_{j}\xi
_{k-j}^{\prime}+a_{k}\xi_{0}^{\prime}+
\sum_{j>k}a_{j}\xi_{k-j} \Biggr)
\Biggr] .
\]
By using definition (\ref{modulusbound}),
\[
\bigl|P_{0}(Y_{k})\bigr|\leq C{\mathbb{E}}_{\xi}\bigl|a_{k}
\bigl(\xi_{0}-\xi_{0}^{\prime
} \bigr)\bigr|^{\gamma}
\bigl(\bigl|X_{k}^{\prime}\bigr|\vee\bigl|X_{k}^{\prime\prime}\bigr|
\bigr)^{\beta},
\]
where $X_{k}^{\prime}=\sum_{j=0}^{k-1}a_{j}\xi_{k-j}^{\prime}+a_{k}\xi
_{0}+\sum_{j>k}a_{j}\xi_{k-j}$ and $X_{k}^{\prime\prime}=\sum_{j=0}^{k-1}a_{j}\xi
_{k-j}^{\prime}+a_{k}\xi_{0}^{\prime}+\sum_{j>k}a_{j}\times\allowbreak\xi_{k-j}$.
Therefore, by
taking the expected value, noticing that $X_{k}^{\prime}$ and $X_{k}^{\prime\prime}$ are
identically distributed as $X_{k}=\sum_{j=0}^{\infty}a_{j}\xi_{k-j}$, and then
applying the Cauchy--Schwarz inequality, for a positive constant~$C^{\prime}$,
we obtain
\[
\bigl\Vert P_{0}(X_{k})\bigr\Vert_{2}^{2}\leq
C^{\prime}a_{k}^{2}\mathbb{E} \bigl(|\xi
_{0}|^{2\gamma} \bigr)\mathbb{E} \bigl(|X_0|^{2\beta}
\bigr).
\]
We estimate now $\mathbb{E}(|X_{0}|^{2\beta})$. If $\beta<1$ then
$\mathbb{E}(|X_{0}|^{2\beta})\leq(\mathbb{E}|X_{0}|^{2})^{\beta}\leq
(\sum_{l\geq 0}a_{l}^{2})^{\beta}\sigma^{2\beta}$.
In case $\beta\geq1$, by the
Rosenthal inequality (see Theorem~1.5.9 in \cite{dePG}), for some positive
constant $C_{\beta}$,  $\mathbb{E}(|X_{0}|^{2\beta})\leq C_{\beta}
((\sum_{l\geq0}a_{l}^{2})^{\beta}\sigma^{2\beta}+\sum_{l\geq0}a_{l}^{2\beta}\mathbb{E}(|\xi_{0}|^{2\beta}))$.
Since we assume that $\sum_{l\geq0}a_{l}^{2}<\infty$, it follows that we can find a constant $K$ such that
\[
\bigl\Vert P_{0}(Y_{k})\bigr\Vert_{2}^{2}\leq
Ka_{k}^{2}\mathbb{E} \bigl(|\xi_{0}|^{2\gamma}
\bigr) \bigl(\mathbb{E}|\xi_{0}|^{2\beta}\vee
\sigma^{2\beta} \bigr).
\]
The result follows by (\ref{condproj}) and by taking into account condition
(\ref{condcorflin}).
\end{pf}

\subsection{Application to mixing stationary sequences}\label{s35}

Mixing coefficients are important for quantifying the strength of dependence
in a stochastic process. They have proven essential for analyzing Markov
chains, Gaussian processes, dynamical systems and other dependent structures.

We shall introduce the following strong mixing coefficient:
For a $\sigma$-algebra $\mathcal{A}$ and a random variable $X$ the strong mixing
coefficient is defined as
\[
\tilde{\alpha}(\mathcal{A},X)=\sup \bigl\{\bigl|\mathbb{P} \bigl(A \cap \{X>x \} \bigr)-\mathbb{P}(A)\mathbb{P}(X>x)\bigr|;x\in R \bigr\}.
\]
This coefficient was introduced by Rosenblatt \cite{Ros56} and also analyzed
by Rio \cite{Rio}. It is weaker than those involving all the future of the
process which are usually used in the literature and they are estimable for a
variety of examples from dynamical systems.

For a stationary sequence of random variables
$(X_{k})_{k\in\mathbb{Z}}$, we
denote by $\mathcal{F}_{m}$ the $\sigma$-field generated by $X_{l}$ with
indices $l\leq m$. Notice that $(\mathcal{F}_{k})_{k\in\mathbb{Z}}$ defined in
this way is a minimal filtration such that $(X_{k})_{k\in\mathbb{Z}}$ is
adapted to $(\mathcal{F}_{k})_{k\in\mathbb{Z}}$. The sequences of coefficients
$\tilde{\alpha}(n)$ are then defined by
\[
\tilde{\alpha}(n)=\tilde{\alpha}(\mathcal{F}_{0},X_{n}).
\]
We refer to the book by Bradley \cite{Bradley} for classical mixing
coefficients and to Dedecker \textit{et al.} \cite{DGM} for specific estimates
of coefficients of type $\tilde{\alpha}$ for certain dynamical systems
generated by intermittent maps.

For integrable random variable $X_{0}$, define the ``upper
tail'' quantile function $Q$ by
\[
Q(u)=\inf \bigl\{ t\geq0\dvtx \mathbb{P} \bigl( |X_{0}|>t \bigr) \leq u \bigr\}.
\]
By relation (1.11c) in Rio \cite{Rio} notice that
\begin{equation}\label{boundalfa}
\bigl\Vert \mathbb{E}_{0}(X_{k})\bigr\Vert_{2}^{2}=\mathbb{E} \bigl(X_{k}{\mathbb{E}}
_{0}(X_{k})
\bigr)\leq 2\int_{0}^{\tilde{\alpha}(k)}Q^{2}(u)\,
\mathrm{d}u.
\end{equation}
By using this inequality, condition (\ref{condcorquenched}) is verified
provided
\begin{equation}\label{strong}
\sum_{k=1}^{\infty}\frac{1}{k}\int
_{0}^{\tilde{\alpha}(k)}Q^{2}(u)\,\mathrm{d}u<\infty.
\end{equation}
Denoting $\tilde{\alpha}^{-1}(x)=\min\{k\in{\mathbb{N}}
 \dvtx  \tilde{\alpha}(k)\leq x\}$ we can write relation (\ref{strong}) in the
equivalent formulation
\[
\int_{0}^{1}\log \bigl(1+\tilde{
\alpha}^{-1}(u) \bigr)Q^{2}(u)\,\mathrm{d}u<\infty.
\]
In particular, if $E(|X_{0}|^{2+\delta})<\infty$ for some positive $\delta>0$,
by decoupling the above integral via the Cauchy--Schwarz inequality, we obtain
that a sufficient condition for (\ref{condcorquenched}) is
\[
\int_{0}^{1} \bigl[\log \bigl(1+\tilde{
\alpha}^{-1}(u) \bigr) \bigr]^{(2+\delta)/\delta}\,\mathrm{d}u<\infty,
\]
which requires a logarithmic rate of decay of the coefficients $\tilde{\alpha
}(k)$. If $\Vert X_{0}\Vert_{\infty}<\infty$, condition (\ref{strong}) is implied by
\[
\sum_{k=1}^{\infty}\frac{1}{k}\tilde{\alpha}(k)<\infty.
\]
Since $(\tilde{\alpha}(k))_{k\geq1}$ is decreasing, by (\ref{boundalfa}),
condition (\ref{strong}) implies the regularity condition (\ref{regular}).
Therefore it is a sufficient condition for (\ref{QCLT1}) which holds with
$N_{1}(t)$ and $N_{2}(t)$ i.i.d. normal random variables, mean $0$
and variance  $\uppi f(t)$, $f(t)$ being the spectral density of the process.

It is worth mentioning that some more restrictive mixing conditions make
possible to obtain (\ref{QCLT1}) directly from (\ref{ACLT}). One of these
conditions is called $\phi$-mixing. A stationary sequence of random variables
$(X_{k})_{k\in\mathbb{Z}}$ is called $\phi$-mixing if
\[
\phi(n)=\sup \bigl\{\bigl|\mathbb{P}(B|A)-\mathbb{P}(B)\bigr|;A\in\mathcal{F}_{0},
B\in\mathcal{F}^{n} \bigr\}\rightarrow0.
\]
Here $\mathcal{F}^{n}$ is the $\sigma$-field generated by $X_{l}$ with
indices $l\geq n$. It is equivalent to saying that (see \cite{Bradley}, Vol.~1)
\[
\phi(n)=\sup \bigl\{\bigl|\mathbb{P}(B|\mathcal{F}_{0})-\mathbb{P}(B)\bigr|;
B\in\mathcal{F}^{n} \bigr\}\rightarrow 0 \qquad\mbox{a.s.}
\]
If we fix now $m>0$,  we have $S_{m}(t)/\sqrt{n}\rightarrow 0$ $\mathbb{P}$-a.s. and it is enough to study the asymptotic behavior of
\[
\mathbf{V}_{n,m}(t)= \bigl(\operatorname{Re} \bigl[S_{n}(t)-S_{m}(t)
\bigr]/\sqrt{n},\operatorname{Im} \bigl[S_{n}(t)-S_{m}(t)
\bigr]/\sqrt{n} \bigr).
\]
By the definition of $\phi$-mixing coefficients, for $h$ continuous and
bounded (see again \cite{Bradley}, Vol.~1)
\[
\bigl|\mathbb{E} \bigl(h \bigl(\mathbf{V}_{n,m}(t) \bigr)|
\mathcal{F}_{0} \bigr)-\mathbb{E} \bigl(h \bigl(\mathbf{V}_{n,m}(t)
\bigr) \bigr)\bigr|\leq\phi(m)\qquad\mbox{a.s.},
\]
and the claim follows easily by Theorem~3.2 in \cite{Bil}.

\section{Proofs}\label{s4}

\begin{pf*}{Proof of Theorem~\ref{c-quenched}}
The proof of Theorem~\ref{c-quenched}\textbf{ }is based on the following
approximation lemma for Fourier transforms. Recall the definition of
projection operator (\ref{projop}).

\begin{lemma}
\label{mart-apprx-Fourier}
Under the conditions of Theorem~\ref{c-quenched},
for almost all $t\in [0,2\uppi]$, the martingale difference
\[
D_{k}(t,\omega)=\sum_{j>k}
\mathrm{e}^{\mathrm{i}jt}P_{k}X_{j}(\omega)=
\mathrm{e}^{\mathrm{i}kt}\sum_{j\geq1}
\mathrm{e}^{\mathrm{i}jt}P_{0}X_{j}(\omega)\circ
T^{k}
\]
is well defined in the almost sure sense and in $\mathbb{L}^{2}(\Omega
,\mathcal{K},\mathbb{P})$. Denote by $M_{n}(t)(\omega)=\sum_{k=1}^{n}D_{k}(t,\allowbreak\omega)$.
Then, for almost all $t\in [0,2\uppi]$,
\[
\frac{1}{n}\mathbb{E}_{0}\bigl|S_{n}(t)-
\mathbb{E}_{0} \bigl(S_{n}(t) \bigr)-M_{n}
(t)\bigr|^{2}
\rightarrow 0\qquad\mathbb{P}\mbox{-a.s. and in }\mathbb{L}_{1}.
\]
\end{lemma}

\begin{pf}
The convergence in $\mathbb{L}_{1}$ was established in
Peligrad and Wu \cite{PeWu}. We shall prove here the almost sure convergence.
It is convenient to work on the product space, $(\widetilde{\Omega},
\widetilde{\mathcal{F}},\tilde{\mathbb{P}})=([0,2\uppi]\times
\Omega,{\mathcal{B}}\otimes{\mathcal{A}},\lambda\otimes{\mathbb{P}})$ where
$\lambda$ is the normalized Lebesgue measure on $[0,2\uppi]$, and ${\mathcal{B}}$
is the Borel $\sigma$-algebra on $[0,2\uppi]$, $\tilde{\mathbb{P}}=\lambda\otimes{\mathbb{P}}$.
Consider also\vspace*{1pt} the filtration $(\widetilde{\mathcal{F}}_{n})_{n\in\mathbb{Z}}$ given by $\widetilde{\mathcal{F}}_{n}:=
\mathcal{B}\otimes {\mathcal{F}}_{n}$. Denote by\vspace*{1pt} $\tilde{\mathbb{E}}$,
the integral with respect to $\tilde{\mathbb{P}}$,  by $\tilde{\mathbb{E}}_{0}$
the conditional expectation with respect to $\widetilde{\mathcal{F}}_{0}$,
$\tilde{P}_{k}(\cdot)=\tilde{\mathbb{E}}_{k}(\cdot)-\tilde{\mathbb{E}}_{k-1}(\cdot)$.

Let $t\in [0,2\uppi)$ be a real number, fixed for the moment. Clearly, the
transformation $\widetilde{T}_{t}$ from $\widetilde{\Omega}$ to $\widetilde
{\Omega}$ given by
\[
\widetilde{T}_{t} \dvtx  (u,\omega)\mapsto \bigl(u+t\,\mbox{modulo}\,2
\uppi,T(\omega) \bigr),
\]
is invertible, bi-measurable and preserves $\tilde{\mathbb{P}}$.
For every
$(u,\omega)\in\widetilde{\Omega}$ define the variable $\tilde{X}_{0}$
on~$\widetilde{\Omega}$ by $\tilde{X}_{0}(u,\omega)=\mathrm{e}^{\mathrm{i}u}
X_{0}(\omega)$ and for any $n\in{\mathbb{Z}}$, $\tilde{X}_{n}(t;u,\omega)
=\tilde{X}_{0}(u,\omega)\circ\widetilde{T}_{t}^{n}$. For
simplicity, in the sequel, we shall drop from the notation the variables $u$
and $\omega$ in $\tilde{X}_{k}(t;u,\omega)$ and we shall write instead
$\tilde{X}_{k}(t)$ and $\tilde{S}_{n}(t)=\sum_{k=1}^{n}\tilde{X}_{k}(t)$.
Notice that $(\tilde{X}_{n}(t))_{n\in{\mathbb{Z}}}$ is a
stationary sequence of complex random variables adapted to the nondecreasing
filtration $(\widetilde{\mathcal{F}}_{n})$.\vspace*{1pt} Moreover $\mathrm{e}^{\mathrm{i}u}\mathrm{e}^{\mathrm{i}kt}X_{k}(\omega)=\tilde{X}_{k}(t;u,\omega)$. We shall
construct a martingale $\tilde{M}_{n}(t)$,
adapted to $(\widetilde{\mathcal{F}}_{n})$,
with stationary differences, such that for almost all
$t\in [0,2\uppi]$
\[
\frac{1}{n}\tilde{\mathbb{E}}_{0} \bigl[\tilde{S}_{n}(t)-
\tilde{\mathbb{E}}_{0} \bigl(\tilde{S}_{n}(t)
\bigr)-\tilde{M}_{n}(t) \bigr]^{2}\rightarrow0\qquad\tilde{\mathbb{P}}\mbox{-a.s.}
\]
With this aim we shall apply Proposition~\ref{mart-approx}, given in the
Section~\ref{s5}. In order to verify the conditions of this proposition, we have to
show that for almost all $t$ in $[0,2\uppi]$
\begin{equation}\label{cond2}
\tilde{P}_{0} \bigl(\tilde{S}_{n}(t) \bigr)\rightarrow
\tilde{D}_{0}(t)\qquad\tilde{\mathbb{P}}\mbox{-a.s.}
\end{equation}
and
\begin{equation}\label{cond3}
\tilde{\mathbb{E}} \Bigl[\sup_{n}\bigl|\tilde{P}_{0}
\bigl(\tilde{S}_{n}(t) \bigr)\bigr|^{2} \Bigr]<\infty.
\end{equation}
In order to prove (\ref{cond2}), note that by the orthogonality of the
projections and the fact that the sequence $\Vert E_{-n}X_{0}\Vert_{2}$ is
decreasing, it follows that
\begin{eqnarray*}
\sum_{k\geq0}\Vert P_{-k}X_{0}
\Vert_{2}^{2} & =& \lim_{n\rightarrow\infty}
\sum_{k=0}^{n}\Vert P_{-k}X_{0}
\Vert_{2}^{2}=\lim_{n\rightarrow\infty}
\biggl\Vert\sum_{k=0}^{n}P_{-k}X_{0}
\biggr\Vert_{2}^{2}
\\
& =& \lim_{n\rightarrow\infty}\Vert X_{0}-E_{-n}X_{0}
\Vert_{2}^{2}=\Vert X_{0}\Vert_{2}^{2}-
\Vert E_{-\infty}X_{0}\Vert_{2}^{2}\leq\Vert
X_{0}\Vert _{2}^{2}<\infty.
\end{eqnarray*}
Clearly this implies
\[
\sum_{k\geq0}|P_{0}X_{k}|^{2}<
\infty\qquad \mathbb{P}\mbox{-a.s.}
\]
Now\vspace*{1pt} for $\omega$ such that $\sum_{k\geq0}|P_{0}X_{k}|^{2}(\omega)<\infty$,
by Carleson's \cite{Ca} theorem, $P_{0}S_{n}(t)=  \sum_{1\leq k\leq n}%
\mathrm{e}^{\mathrm{i}kt}(P_{0}X_{k})(\omega)$ converges $\lambda$-almost surely.
Denote the limit by $D_{0}=D_{0}(t)$. We now consider the set
\[
A= \bigl\{(t,\omega)\subset [0,2\uppi]\times\Omega,\mbox{ where }
\bigl[P_{0}S_{n}(t) \bigr]_{n}\mbox{ does not converge} \bigr\}
\]
and note that almost all sections for $\omega$ fixed have Lebesgue measure
$0$. So by Fubini's theorem the set $A$ has measure $0$ in the product space
and therefore, again by Fubini's theorem, almost all sections for $t$ fixed
have probability $0$. It follows that for almost all $t$ in $[0,2\uppi]$,
$P_{0}(S_{n}(t))\rightarrow D_{0}(t)$ almost surely under $\mathbb{P}$. This
shows that, after multiplying by $\mathrm{e}^{\mathrm{i}u}$, we get, for
almost all $t$, that condition (\ref{cond2}) is verified with
\[
\tilde{D}_{0}(t)=\mathrm{e}^{\mathrm{i}u}\sum
_{j\geq1}\mathrm{e}^{\mathrm{i}jt}P_{0}X_{j}(\omega).
\]
Note that
\[
\tilde{D}_{k}(t)=\tilde{D}_{0}(t)\circ
\widetilde{T}_{t}^{k}
=\mathrm{e}^{\mathrm{i}u}
\sum_{j>k}\mathrm{e}^{\mathrm{i}jt}P_{k}X_{j}(
\omega).
\]
Next, we prove (\ref{cond3}). By the maximal inequality in Hunt and Young
\cite{HY}, there is a constant $C$ such that
\[
\int_{0}^{2\uppi} \Bigl[ \sup_{n}\bigl|P_{0}
\bigl(S_{n}(t) \bigr)\bigr|^{2} \Bigr] \lambda(\mathrm{d}t)\leq C
\sum_{k\geq1}|P_{0}X_{k}|_{2}^{2}.
\]
Then we integrate with respect to $\mathbb{P}$ and use Fubini theorem to
obtain
\[
\int_{0}^{2\uppi}\mathbb{E} \Bigl[ \sup
_{n\geq0}\bigl|P_{0} \bigl(S_{n}(t)
\bigr)\bigr|^{2} \Bigr] \lambda(\mathrm{d}t)\leq C\Vert X_{0}
\Vert_{2}^{2}<\infty.
\]
It follows that
\[
\mathbb{E} \Bigl[ \sup_{n\geq1}\bigl|P_{0}
\bigl(S_{n}(t) \bigr)\bigr|^{2} \Bigr] <\infty\qquad\mbox{for almost all }t.
\]
Therefore, we obtain that condition (\ref{cond3}) is satisfied. We apply now
Proposition~\ref{mart-approx} to obtain for almost all $t$ in $[0,2\uppi]$
\begin{equation}\label{martbarapprox}
\frac{1}{n}\tilde{\mathbb{E}}_{0}\bigl|\tilde{S}_{n}(t)-
\tilde{\mathbb{E}}
_{0} \bigl(\tilde{S}_{n}(t)
\bigr)-\tilde{M}_{n}(t)\bigr|^{2}\rightarrow 0 \qquad
\tilde{\mathbb{P}}\mbox{-a.s.},
\end{equation}
where
\[
\tilde{M}_{n}(t)=\sum_{k=1}^{n}
\tilde{D}_{k}(t).
\]
Now fix $t$ in $[0,2\uppi]$ such that (\ref{martbarapprox}) holds. Clearly
\[
\frac{1}{n}\mathbb{E}_{0}\bigl|S_{n}(t)-
\mathbb{E}_{0} \bigl(S_{n}(t) \bigr)-M_{n}%
(t)\bigr|^{2}
\rightarrow 0\qquad \mathbb{P}\mbox{-a.s.}
\]
The result follows.
\end{pf}

We study next the behavior of $M_{n}(t)/\sqrt{n}$. We shall do it in general
in the context of stationary and ergodic complex valued martingale
differences. Below, the martingale difference $D$ may depend on $t$.

\begin{proposition}
\label{Q-martingale}
Let $T$ and $\mathcal{F}_{0}$ be as in Section~\ref{s2}. Assume
that $t\in(0,2\uppi)$ be such that $\mathrm{e}^{-2\mathrm{i}t}$ is not an eigenvalue of
$T$. Let $D=D_{0}$ be a random variable defined on $(\Omega,\mathcal{F},\mathbb{P})$,
$\mathcal{F}_{0}$ measurable and such that $E(D\circ
T|\mathcal{F}_{0})=0$ a.s. For any $k\in{\mathbb{Z}}$, let
\[
\mathbf{D}_{k}(t)= \bigl(\operatorname{Re} \bigl(
\mathrm{e}^{\mathrm{i}kt}D \circ T^{k} \bigr),
\operatorname{Im} \bigl( \mathrm{e}^{\mathrm{i}kt}D\circ T^{k} \bigr)
\bigr).
\]
Let $\mathbf{M}_{n}(t)=\sum_{k=1}^{n}\mathbf{D}_{k}(t)$. Then, there is a set
$\Omega^{\prime}$ with $\mathbb{P}(\Omega^{\prime})=1$ such that for all
$\omega\in\Omega^{\prime}$
\begin{equation}\label{CLTmart}
\frac{1}{\sqrt{n}}\mathbf{M}_{n}(t)\Rightarrow\mathbf{N}\qquad\mbox{under }
\mathbb{P}^{\omega},
\end{equation}
where $\mathbf{N}=(N_{1},N_{2})$, with $N_{1}, N_{2}$ are two
independent centered normal random variables with variance $\mathbb{E}|D|^{2}/2$.
\end{proposition}

\begin{pf}
Fix $t\in(0,2\uppi)$ such that $\mathrm{e}^{-2\mathrm{i}t}$ is not an
eigenvalue of $T$. Denote $R_{k}(t)=\operatorname{Re}(\mathrm{e}^{\mathrm{i}kt}D\circ
T^{k})$ and $I_{k}(t)=\operatorname{Im}(\mathrm{e}^{\mathrm{i}kt}D\circ T^{k})$.

The proof is based on Theorem~\ref{martCLT} and the following two convergence
results: for any real constants $a$ and $b$
\begin{equation}\label{1}
\frac{1}{\sqrt{n}}\mathbb{E}_{0} \Bigl(\max_{1\leq k\leq n}\bigl|aR_{k}(t)+bI_{k}(t)\bigr|
\Bigr)\rightarrow 0\qquad \mathbb{P}\mbox{-a.s.}
\end{equation}
and
\begin{equation}\label{2}
\mathbb{P}_{0} \Biggl(\Biggl|\frac{1}{n}\sum
_{k=1}^{n}\bigl|aR_{k}(t)+bI_{k}(t)\bigr|^{2}-
\frac{1}{2} \bigl(a^{2}+b^{2} \bigr)
\mathbb{E}|D|^{2}\biggr|> \varepsilon \Biggr)\rightarrow 0 \qquad \mathbb{P}
\mbox{-a.s.}
\end{equation}
Before proving (\ref{1}) and (\ref{2}) let us show how they lead to the result.

Let $a$ and $b$ be two rational numbers and let $\Omega_{a,b}$ be the set of
probability $1$ where (\ref{1}) and~(\ref{2}) hold. Construct $\Omega_{1}=\bigcap\Omega_{a,b}$,
where the intersection is taken over all the rationals $a$
and $b$. Clearly $\mathbb{P}(\Omega_{1})=1$. Then, by Theorem~\ref{martCLT} in
Section~\ref{s5}, we get via (\ref{1}) and (\ref{2}) that for all $\omega\in
\Omega_{1}$
\begin{equation}\label{CW}
\sum_{k=1}^{n} \bigl(aR_{k}(t)+bI_{k}(t)
\bigr)\Big/\sqrt{n}\Rightarrow N(a,b,t)\qquad\mbox{under }\mathbb{P}^{\omega},
\end{equation}
where $N(a,b,t)$ is a centered normal random variable with variance
$(a^{2}+b^{2})\mathbb{E}|D|^{2}/2$.

Because $\mathbb{E}_{0}$ is regular, by Hopf ergodic theorem
\[
\frac{1}{n}\mathbb{E}_{0}\bigl|\mathbf{M}_{n}(t)\bigr|^{2}=
\frac{1}{n}\sum_{k=1}^{n}
\mathbb{E}_{0}\bigl|\mathbf{D}_{k}(t)\bigr|^{2}\rightarrow
\mathbb{E}\bigl|\mathbf{D}_{0}(t)\bigr|^{2}\qquad\mbox{as }n
\rightarrow\infty  \ \mathbb{P}\mbox{-a.s.}
\]
By Markov inequality it follows that there is a set $\Omega_{2}$ such that for
all $\omega\in\Omega_{2}$ the sequence $(\mathbf{M}_{n}(t)/\sqrt{n})_{n\geq1}$
is tight under $\mathbb{P}^{\omega}$.

Now construct $\Omega^{\prime}=\Omega_{1}\cap\Omega_{2}$. For $\omega\in
\Omega^{\prime}$, we apply Lemma~\ref{weak-conv} in Section~\ref{s5} and obtain
(\ref{CLTmart}).

It remains to prove (\ref{1}) and (\ref{2}). To prove the convergence in
(\ref{2}) we shall use relation (16) in Cuny--Merlev\`{e}de--Peligrad
\cite{CMP}, with $u=0$, which gives
\[
\frac{1}{n}\sum_{k=1}^{n}\bigl|aR_{k}(t)+bI_{k}(t)\bigr|^{2}
\rightarrow\frac{1}{2}
 \bigl(a^{2}+b^{2}
\bigr)\mathbb{E}|D|^{2}\qquad \mathbb{P}\mbox{-a.s.}
\]
This convergence was obtained by trigonometric computations along with
Dunford--Schwartz ergodic theorem from Sections~VIII.5 and~VIII.6 of \cite{DS},
which requires that $t\in(0,2\uppi)$ be such that $\mathrm{e}^{-2\mathrm{i}t}$ is not an
eigenvalue of $T$ (see Proposition~30 in \cite{CMP}).

This last convergence implies that, for every $\varepsilon>0$
\[
I \Biggl(\Biggl|\frac{1}{n}\sum_{k=1}^{n}\bigl|aR_{k}(t)+bI_{k}(t)\bigr|^{2}-
\frac{1}{2} \bigl(a^{2}
+b^{2} \bigr)
\mathbb{E}|D|^{2}\Biggr|>\varepsilon \Biggr)\rightarrow 0 \qquad \mathbb{P}\mbox{-a.s.},
\]
whence (\ref{2}) follows by Theorem~34.2(v) in Billingsley \cite{BilPrM}.

We verify now relation (\ref{1}). Note that
\[
\bigl|aR_{k}(t)+bI_{k}(t)\bigr|\leq\bigl(|a|+|b|\bigr)|D|\circ
T^{k}=\bigl(|a|+|b|\bigr)|D_{k}|.
\]
It is enough to verify that
\[
\frac{1}{n}\mathbb{E}_{0} \Bigl(\max_{1\leq k\leq n}|D_{k}|^{2}
\Bigr)\rightarrow 0 \qquad \mathbb{P}\mbox{-a.s.}
\]
We shall use a truncation argument. Let $\varepsilon>0$ and $c>0$ be fixed for
the moment. Let $n$ be sufficiently large such that $\varepsilon\sqrt{n}\geq
c$.  For this selection of $n$, we have
\begin{eqnarray*}
\frac{1}{n}\mathbb{E}_{0} \Bigl(\max_{1\leq k\leq n}|D_{k}|^{2}
\Bigr) & \leq & \frac{1}{n}\mathbb{E}_{0} \Bigl(\max
_{1\leq k\leq n}|D_{k}|^{2}I\bigl(|D_{k}|\leq
\varepsilon \sqrt{n}\bigr) \Bigr)
\\
&&{}+\frac{1}{n}\mathbb{E}_{0} \Bigl(\max_{1\leq k\leq n}|D_{k}|^{2}I\bigl(|D_{k}|
>\varepsilon\sqrt{n}\bigr) \Bigr)\\
&\leq& \varepsilon^{2}+\frac{1}{n}
\sum_{k=1}^{n}
\mathbb{E}_{0} \bigl(|D_{k}|^{2}I\bigl(|D_{k}|>c\bigr)
\bigr).
\end{eqnarray*}
Now, by the Hopf theorem for Dunford--Schwartz operators (see \cite{DS} or
\cite{E}),
\[
\frac{1}{n}\sum_{k=1}^{n}
\mathbb{E}_{0} \bigl(|D_{k}|^{2}I\bigl(|D_{k}|>c\bigr)
\bigr)\rightarrow \mathbb{E} \bigl(|D_{0}|^{2}I\bigl(|D_{0}|>c\bigr)
\bigr)\qquad \mbox{as }n\rightarrow\infty \ \mathbb{P}\mbox{-a.s.}
\]
Then we have
\[
\lim\sup_{n\rightarrow\infty}\frac{1}{n}\mathbb{E}_{0}
\Bigl(\max_{1\leq k\leq
n}|D_{k}|^{2} \Bigr)\leq
\varepsilon^{2}+\mathbb{E} \bigl(|D_{0}|^{2}I\bigl(|D_{0}|>c\bigr)
\bigr).
\]
The result follows by letting $\varepsilon\rightarrow0$ and $c\rightarrow
\infty$.
\end{pf}

\begin{remark}
\label{remQmart}
Note that because $\mathcal{K}$ is countably generated then
$\mathbb{L}^{2}(\Omega,\mathcal{K},\mathbb{P})$ is separable and by Lemma~32
in Cuny \textit{et al.} \cite{CMP}, $T$ can admit only a countable number of
eigenvalues. Therefore the quenched CLT in Proposition~\ref{Q-martingale}
holds for almost all $t\in [0,2\uppi]$.
\end{remark}

\begin{pf*}{End of the Proof of Theorem~\ref{c-quenched}}
By using Theorem~3.1 in Billingsley \cite{Bil}, Lemma~\ref{mart-apprx-Fourier}
shows that for almost all $t\in [0,2\uppi]$, there is a set $\Omega
^{\prime}\subset\Omega$ with $\mathbb{P}(\Omega^{\prime})=1$ such that for all
$\omega\in\Omega^{\prime}$, the limiting behavior $S_{n}(t)-\mathbb{E}_{0}(S_{n}(t))$ is the same as of the martingale $M_{n}(t)$ under
$\mathbb{P}^{\omega}$. Then, by Proposition~\ref{Q-martingale} and Remark~\ref{remQmart}, for almost all $t\in [0,2\uppi]$ the quenched CLT holds for
$\mathbf{W}_{n}(t)$, with the limit $\mathbf{N}(t)=(N_{1}(t), N_{2}(t))$,
where $N_{1}(t), N_{2}(t)$ are two independent centered normal random
variables with variance $\mathbb{E}|D(t)|^{2}/2$. For an alternative
characterization of $\mathbb{E}|D(t)|^{2}$, it remains to note that by Lemma~\ref{mart-apprx-Fourier}, for almost all $t\in [0,2\uppi]$
\[
\biggl(\frac{1}{n}\mathbb{E}_{0}\bigl|S_{n}(t)-
\mathbb{E}_{0} \bigl(S_{n}(t) \bigr)\bigr|^{2}
 \biggr)^{1/2}- \biggl(\frac{1}{n}\mathbb{E}_{0}\bigl|M_{n}(t)\bigr|^{2}
\biggr)^{1/2}\rightarrow0\qquad \mathbb{P}\mbox{-a.s.}
\]
Furthermore, by the Hopf ergodic theorem
\[
\frac{1}{n}\mathbb{E}_{0}\bigl|M_{n}(t)\bigr|^{2}=
\frac{1}{n}\sum_{k=1}^{n}%
\mathbb{E}_{0}\bigl|D_{k}(t)\bigr|^{2}\rightarrow
\mathbb{E}\bigl|D(t)\bigr|^{2} \qquad
\mathbb{P}\mbox{-a.s.}
\]
\upqed\end{pf*}
\noqed\end{pf*}

\begin{pf*}{Proof of Theorem~\ref{char}}
Clearly (b) implies (a) via Theorem~\ref{c-quenched}. To prove that (a)
implies (b), we shall use again Theorem~\ref{c-quenched} along with the Theorem
of types. This latter theorem states that if $V_{n}=a_{n}U_{n}+b_{n}$ and
$V_{n}\Rightarrow V$ and $U_{n}\Rightarrow U$ with $U$ nondegenerate then
$a_{n}\rightarrow a$, $b_{n}\rightarrow b$ and $V=aU+b$.

Under conditions of Theorem~\ref{char}, for $\lambda$-almost all $t\in
[0,2\uppi]$ there is a set $\Omega^{\prime}\subset\Omega$
with
\mbox{$\mathbb{P}(\Omega^{\prime})=1$} such that for all $\omega\in\Omega^{\prime}$
\[
\frac{1}{\sqrt{n}}\operatorname{Re} \bigl[S_{n}(t)-
\mathbb{E}_{0} \bigl(S_{n}(t) \bigr) \bigr]
\Rightarrow N_{1}(t)\qquad \mbox{under }\mathbb{P}^{\omega}.
\]
By the properties of conditional expectations and measure theoretical
arguments (see Lemma~\ref{reg}), we know that for every function $g$
continuous and bounded and random variables $X$ and $Y$, such that $Y$ is
$\mathcal{F}_{0}$-measurable,
\[
\mathbb{E}^{\omega} \bigl(g(X,Y)|\mathcal{F}_{0} \bigr)=
\mathbb{E}^{\omega} \bigl(g \bigl(X,Y(\omega ) \bigr)|\mathcal{F}_{0}
\bigr)
\]
for $\omega$ in a set of probability $1$. By this observation along with the
definition of convergence in distribution, we derive that for $\lambda$-almost
all $t\in [0,2\uppi]$ there is a set $\Omega^{\prime}\subset\Omega$ with
$\mathbb{P}(\Omega^{\prime})=1$, such that for all $\omega\in\Omega^{\prime}$
\[
\frac{1}{\sqrt{n}}\operatorname{Re} \bigl[S_{n}(t)-
\mathbb{E}_{0} \bigl(S_{n}
(t) \bigr) (\omega)
\bigr]\Rightarrow N_{1}(t)\qquad\mbox{under }\mathbb{P}^{\omega},
\]
and by (a) there is a set $\Omega^{\prime\prime}\subset\Omega$ with $\mathbb{P}(\Omega^{\prime\prime})=1$,
such that for all $\omega\in\Omega^{\prime\prime}$
\[
\frac{1}{\sqrt{n}}\operatorname{Re} \bigl(S_{n}(t) \bigr)\Rightarrow
N_{1}(t)\qquad\mbox{under }\mathbb{P}^{\omega}.
\]
Now assume that $N_{1}(t)$ is nondegenerate. For $\omega\in\Omega^{\prime}\cap\Omega^{\prime\prime}$,
by the Theorem of types we have $\operatorname{Re}\mathbb{E}_{0}(S_{n}(t))(\omega)/\sqrt{n}\rightarrow0$. A similar argument
gives $\operatorname{Im}\mathbb{E}_{0}(S_{n}(t))(\omega)/\sqrt{n}\rightarrow0$
and (b) follows for this case. If $N_{1}(t)$ is degenerate, then both
$\operatorname{Re}[S_{n}(t)]/\sqrt{n}\rightarrow0$ under $\mathbb{P}^{\omega}$
and $\operatorname{Re}[S_{n}(t)-\mathbb{E}_{0}(S_{n}(t))(\omega)]/\sqrt
{n}\rightarrow0$ under $\mathbb{P}^{\omega}$, and the result follows.
\end{pf*}

\begin{pf*}{Proof of Corollary~\ref{corquenched}}
In order to prove this result,
we shall verify the item (b) of Theorem~\ref{char}. We have then to show that
for almost all $t$
\begin{equation}\label{toshow}
\frac{\mathbb{E}_{0}(S_{n}(t))}{\sqrt{n}}\rightarrow0 \qquad\mathbb{P}\mbox{-a.s.}
\end{equation}
Note that it is enough to show instead that for almost all $t\in [0,2\uppi]$
\[
\frac{(1-\mathrm{e}^{\mathrm{i}t})\mathbb{E}_{0}(S_{n}(t))}{\sqrt{n}}\rightarrow 0\qquad  \mathbb{P}\mbox{-a.s.}
\]
With this aim note that
\begin{eqnarray*}
\frac{(1-\mathrm{e}^{\mathrm{i}t})\mathbb{E}_{0}(S_{n}(t))}{\sqrt{n}} &=& \frac{\mathbb{E}_{0}
(S_{n}(t))-\mathrm{e}^{\mathrm{i}t}\mathbb{E}_{0}(S_{n}(t))}{\sqrt{n}}
\\
&=& \frac{1}{\sqrt{n}}\mathrm{e}^{\mathrm{i}t}\mathbb{E}_{0}(X_{1})-\mathrm{e}^{\mathrm{i}t(n+1)}
\frac{1}{\sqrt{n}}\mathbb{E}_{0}(X_{n})
+\frac{1}{\sqrt{n}}\sum_{k=1}^{n-1}\mathrm{e}^{\mathrm{i}t(k+1)}
\mathbb{E}_{0}(X_{k+1}-X_{k}).
\end{eqnarray*}
We shall analyze each term in the last sum separately. The first term,
$\mathrm{e}^{\mathrm{i}t}\mathbb{E}_{0}(X_{1})/\sqrt{n}$ in the above expression is trivially
convergent to $0$, $\mathbb{P}$-a.s. By Jensen's inequality the
second one is dominated as follows:
\[
\biggl|\mathrm{e}^{\mathrm{i}t(n+1)}\frac{1}{\sqrt{n}}\mathbb{E}_{0}(X_{n})\biggr|^{2}
\leq\frac{1}{n}\mathbb{E}_{0} \bigl(X_{n}^{2}
\bigr).
\]
We write
\[
\frac{1}{n}\mathbb{E}_{0} \bigl(X_{n}^{2}
\bigr)=\frac{1}{n}\sum_{j=1}^{n}
\mathbb{E}_{0} \bigl(X_{j}^{2} \bigr)-
\frac{1}{n}\sum_{j=1}^{n-1}
\mathbb{E}_{0} \bigl(X_{j}^{2} \bigr),
\]
which convergence to $0$, $\mathbb{P}\mbox{-a.s.}$, by the Hopf ergodic
theorem for Dunford--Schwartz operators (see again \cite{E}).

To prove the convergence of the third term, since we assumed
(\ref{condcorquenched}), it follows that
\[
\sum_{k\geq1}\frac{|\mathbb{E}_{0}(X_{k+1}-X_{k})|^{2}}{k}<\infty\qquad \mathbb{P}\mbox{-a.s.}
\]
By Carleson theorem (see \cite{Ca}) it follows that for almost all $t$
\[
\sum_{k\geq1}\frac{\mathrm{e}^{\mathrm{i}tk}\mathbb{E}_{0}(X_{k+1}-X_{k})}{k^{1/2}}\mbox{ converges }\mathbb{P}\mbox{-a.s.}
\]
which implies by Kronecker lemma
\[
\frac{1}{\sqrt{n}}\sum_{k=1}^{n-1}\mathrm{e}^{\mathrm{i}t(k+1)}
\mathbb{E}_{0}(X_{k+1}
-X_{k})
\rightarrow 0 \qquad\mathbb{P}\mbox{-a.s.}
\]
which completes the proof of (\ref{toshow}) and of this corollary.
\end{pf*}

\begin{pf*}{Proof of Theorem~\ref{MW}}
With the notations from the proof of Lemma~\ref{mart-apprx-Fourier},
we note that under condition (\ref{condMW}) we also have
\begin{equation}\label{condMWtilda}
\sum_{k\geq1}\frac{1}{k^{3/2}} \bigl( \tilde{\mathbb{E}}\bigl|\tilde{\mathbb{E}}_{0}
\bigl(\tilde{S}_{k}(t) \bigr)\bigr|^{2} \bigr) ^{1/2}<\infty.
\end{equation}
Then, we can apply directly the martingale approximation in Theorem~2.7 in
Cuny and Merlev\`{e}de \cite{CM} which also remains valid for complex valued
variables. It follows that
\[
\frac{1}{\sqrt{n}}\tilde{\mathbb{E}}_{0}\bigl|\tilde{S}_{n}(t)-
\tilde{M}_{n}(t)\bigr|\rightarrow 0 \qquad  \tilde{\mathbb{P}}\mbox{-a.s. and in }\tilde{\mathbb{L}}^{2},
\]
where $\tilde{M}_{n}$ has stationary complex martingale differences defined
by
\[
\tilde{D}_{j}(t)=\sum_{n\geq0}\sum_{k\geq n}\frac{\tilde{P}_{0}(\tilde{X}_{k}(t))\circ T^{j}}{k+1}.
\]
Whence we obtain
\begin{equation}\label{martapprox2}
\frac{1}{\sqrt{n}}\mathbb{E}_{0}\bigl| \bigl(S_{n}(t)
\bigr)-M_{n}(t))\bigr|\rightarrow0 \qquad \mathbb{P}\mbox{-a.s. and in }\mathbb{L}^{2},
\end{equation}
where the differences of the martingale $M_{n}(t)$ are
\[
D_{j}(t)=\sum_{n\geq0}\sum
_{k\geq n}\frac{\mathrm{e}^{\mathrm{i}tk}P_{0}(X_{k})\circ T^{j}}{k+1}.
\]
It follows that for $\mathbb{P}$-almost all $\omega$, under $\mathbb{P}^{\omega}$, the behavior of $S_{n}(t)/\sqrt{n}$ is equivalent to
$M_{n}(t)/\sqrt{n}$. We have just to apply Proposition~\ref{Q-martingale} to
obtain the quenched CLT, where the limiting independent normal variables have
the variance $\Vert D_{0}(t)\Vert_{2}^{2}/2$. It remains to note that by
(\ref{martapprox2}) we can identify $\Vert D_{0}(t)\Vert_{2}^{2}$ as
\[
\lim_{n\rightarrow\infty}\frac{1}{n}\mathbb{E}\bigl|S_{n}(t)\bigr|^{2}=
\bigl\Vert D_{0}%
(t)\bigr\Vert_{2}^{2}.
\]
\upqed\end{pf*}

\section{Technical results}\label{s5}

First, we prove the following martingale approximation for complex valued
random variables. It is similar to Proposition~7 in Cuny and Peligrad
\cite{CP} but we do not assume ergodicity and the variables are complex
valued. Because there are various changes in the proof we give the proof for completeness.

\begin{proposition}
\label{mart-approx}
Assume that $(X_{k})_{k\in\mathbb{Z}}$ is a stationary
sequence of complex valued random variables and let $S_{n}=\sum_{k=1}^{n}X_{k}$.
Assume
\begin{equation}\label{AS}
P_{0}(S_{n})\rightarrow D_{0} \qquad\mbox{converges
a.s.} \quad \mbox{and}\quad  \mathbb{E} \Bigl[\sup_{m}\bigl|P_{0}(S_{m})\bigr|^{2}
\Bigr]<\infty
\end{equation}
(where $P_{0}$ is defined by (\ref{projop})). Then, $D_{0}$ is a martingale
difference and
\[
\frac{1}{n}\mathbb{E}_{0}\bigl(\bigl|S_{n}-
\mathbb{E}_{0}(S_{n})-M_{n}
\bigr|^2 \bigr)
\rightarrow 0\qquad  \mathbb{P}\mbox{-a.s.},
\]
where $M_{n}=\sum_{k=1}^{n}D_{k}$ with $D_{k}=D_{0}\circ T^{k}$.
\end{proposition}

\begin{pf}
Starting from condition (\ref{AS}), we notice that this
condition implies $P_{0}(S_{n})\rightarrow D_{0}$ in $\mathbb{L}^{2}(\mathbb{P})$.
Since $\mathbb{E}_{-1}[P_{0}(S_{n})]=0$ a.s. we
conclude that $\mathbb{E}_{-1}[D_{0}]=0$ a.s. and therefore
$(D_{k})_{k\geq1}$ is a sequence of martingale differences adapted to
$\mathcal{F}_{k}$. We shall approximate $S_{n}$ by $M_{n}+\mathbb{E}_{0}(S_{n})$.
We use now a traditional decomposition of $S_{n}$ in martingale
differences by using the projections on consecutive sigma algebras:
\[
S_{n}-\mathbb{E}_{0}(S_{n})=
\bigl[S_{n}- \mathbb{E}_{n-1}(S_{n}) \bigr]+ \bigl[
\mathbb{E}
_{n-1}(S_{n})- \mathbb{E}_{n-2}(S_{n})
\bigr]+ \cdots + \bigl[\mathbb{E}_{1}(S_{n})-
\mathbb{E}_{0}(S_{n}) \bigr].
\]
So, we have the martingale decomposition
\[
S_{n}-\mathbb{E}_{0}(S_{n})-M_{n}=\sum_{k=1}^{n} \bigl[P_{k}(S_{n}-S_{k-1})-D_{k}
\bigr].
\]
We write now
\[
P_{k}(S_{n}-S_{k-1})-D_{k}=
\bigl[P_{0}(S_{n-k}) \bigr]\circ T^{k}-D_{0}
\circ T^{k},
\]
and so
\[
\sum_{k=1}^{n} \bigl[P_{k}(S_{n}-S_{k-1})-D_{k}
\bigr]=\sum_{k=1}^{n} \bigl(
P_{0}(S_{n-k})-D_{0} \bigr) \circ
T^{k}=\sum_{k=0}^{n-1} \bigl(
P_{0}(S_{k}
)-D_{0} \bigr) \circ
T^{n-k}.
\]
With the notation
\[
P_{0}(S_{k})-D_{0}=G_{k},
\]
we have
\[
S_{n}-\mathbb{E}_{0}(S_{n})-M_{n}=\sum
_{j=0}^{n-1}G_{j}\circ
T^{n-j}=\sum_{j=1}^{n}G_{n-j}
\circ T^{j}.
\]
By the orthogonality of $G_{j}\circ T^{n-j}$ , we have
\[
\mathbb{E}_{0}\bigl|S_{n}-\mathbb{E}_{0}(S_{n})-M_{n}\bigr|^{2}=
\sum_{j=0}^{n-1}
\mathbb{E}_{0} \bigl(|G_{j}|^{2}\circ
T^{n-j} \bigr).
\]
Let $N$ be fixed. For $n$ sufficiently large, we decompose the last sum into a
sum from $1$ to $N$ and one from $N+1$ to $n$. Then
\begin{eqnarray}\label{decompose1}
\mathbb{E}_{0}\bigl|S_{n}-\mathbb{E}_{0}(S_{n})-M_{n}\bigr|^{2}
& =& \sum_{j=0}^{N}
\mathbb{E}_{0} \bigl(|G_{j}|^{2}\circ
T^{n-j} \bigr)+\sum_{j=N+1}^{n-1}
\mathbb{E}_{0} \bigl(|G_{j}|^{2}\circ
T^{n-j} \bigr)
\nonumber\\[-8pt]
\\[-8pt]
& =& A_{n}(N)+B_{n}(N).
\nonumber
\end{eqnarray}
It is then well known that we have for all $j$ fixed
\[
\frac{1}{n-j}\sum_{u=0}^{n-j}
\mathbb{E}_{0} \bigl(|G_{j}|^{2}\circ
T^{u} \bigr)\mbox{ converges as }n\rightarrow\infty\mbox{ almost
surely and in }\mathbb{L}_{1}.
\]
By writing for all $j$ fixed, $0\leq j\leq N$,
\[
\mathbb{E}_{0} \bigl(|G_{j}|^{2}\circ
T^{n-j} \bigr)=\sum_{j=0}^{n-j}
\mathbb{E}_{0}
 \bigl(|G_{j}|^{2}\circ
T^{j} \bigr)-\sum_{j=0}^{n-j-1}
\mathbb{E}_{0} \bigl(|G_{j}|^{2}\circ
T^{j} \bigr),
\]
it follows easily that
%
\begin{equation}\label{convI}
\frac{A_{n}(N)}{n}\rightarrow 0 \qquad\mbox{as }n\rightarrow\infty \
\mathbb{P}\mbox{-a.s. and in }\mathbb{L}_{1}.
\end{equation}
Now we treat $B_{n}(N)$. We bound this term in the following way,
\[
\frac{B_{n}(N)}{n}\leq\frac{1}{n}\sum_{j=1}^{n}
\mathbb{E}_{0} \Bigl[\sup_{m>N}|G_{m}|^{2}
\circ T^{j} \Bigr].
\]
By the Hopf ergodic theorem and the specification in Section (7) in Dedecker
\textit{et al.} \cite{DMP} we have
\begin{eqnarray*}
&&\lim_{n\rightarrow\infty}\frac{1}{n}\sum_{j=1}^{n}
\Bigl[\sup_{m>N} |G_{m}|^{2}
\circ T^{j} \Bigr]
=\mathbb{E} \Bigl[\sup_{m>N}|G_{m}|^{2}| \mathcal{I} \Bigr]
\\
&& \quad  \mathbb{P}\mbox{-a.s. and in }L_{1},
\end{eqnarray*}
where $\mathcal{I}$ is the invariant sigma field. Since by (\ref{AS})
$\sup_{m>N}|G_{m}|^{2}\rightarrow0$ a.s. as $N\rightarrow\infty$ and
$\sup_{m>N}|G_{m}|^{2}\leq\sup_{m}|G_{m}|^{2}\in\mathbb{L}_{1}$, by
Billingsley \cite{Bil}, Theorem~34.2(v) we also have
\[
\lim_{N\rightarrow\infty}\mathbb{E} \Bigl[\sup_{m>N}|G_{m}|^{2}|
\mathcal{I} \Bigr]=0, \qquad  \mathbb{P}\mbox{-a.s. and in }L_{1},
\]
and therefore
\[
\lim_{N\rightarrow\infty}\lim_{n\rightarrow\infty}\frac{1}{n}\sum
_{j=1}^{n} \Bigl[\sup_{m>N}|G_{m}|^{2}
\circ T^{j} \Bigr]=0, \qquad
\mathbb{P}\mbox{-a.s. and in }\mathbb{L}_{1}.
\]
It follows that
\begin{equation}\label{conv2}
\lim_{N\rightarrow\infty}\lim\sup_{n\rightarrow\infty}
\frac{B_{n}(N)}{n}=0, \qquad \mathbb{P}\mbox{-a.s. and in }\mathbb{L}_{1}.
\end{equation}
\end{pf}

We give below a well-known Raikov type central limit theorem for nonstationary martingales.

The following theorem is a variant of Theorem~3.2 in Hall and Heyde \cite{HH}
(see also G\"{a}nssler and H\"{a}usler \cite{GH}).

\begin{theorem}
\label{martCLT}
Assume $(D_{n,i})_{1\leq i\leq n}$ is an array of square
integrable martingale differences adapted to an array $(\mathcal{F}_{n,i})_{1\leq i\leq n}$ of nested sigma fields. Suppose
\begin{equation}\label{negl1}
\mathbb{E} \Bigl(\max_{1\leq j\leq n}|D_{n,j}| \Bigr)
\rightarrow0 \qquad\mbox{as }
n\rightarrow\infty
\end{equation}
and\vspace*{-6pt}
\begin{equation}\label{conv-2}
\sum_{j=1}^{n}D_{n,j}^{2}
\rightarrow^{\mathbb{P}}\sigma^{2} \qquad \mbox{as }n\rightarrow\infty.
\end{equation}
Then $S_{n}=\sum_{j=1}^{n}D_{n,j}$ converges in distribution to a centered
normal variable with variance $\sigma^{2}$.
\end{theorem}

We give now a result on weak convergence needed for the proof of Proposition~\ref{Q-martingale}.

\begin{lemma}
\label{weak-conv}
Assume that the sequence of random variables $(Y_{n},Z_{n})_{n\geq1}$ is tight and for every rational numbers $a$ and $b$ we have
\[
aY_{n}+bZ_{n}\Rightarrow aN_{1}+bN_{2}.
\]
Then $(Y_{n},Z_{n})\Rightarrow(N_{1},N_{2})$.
\end{lemma}

\begin{pf}
Because $(Y_{n},Z_{n})$ is tight, from any subsequence
$(n^{\prime})$ we can extract another subsequence $(n^{\prime\prime})$ convergent in
distribution to $(L_{1},L_{2})$ say.
By the
Cram\'{e}r--Wold device, it follows that for all real numbers $a$ and $b$ we
have
\[
aY_{n}+bZ_{n}\Rightarrow aL_{1}+bL_{2}.
\]
Therefore for all rational numbers $a$ and $b$ we have
\begin{equation}\label{ab}
\mathbb{E}\mathrm{e}^{\mathrm{i}(aN_{1}+bN_{2})}=\mathbb{E}\mathrm{e}^{\mathrm{i}(aL_{1}+bL_{2})}.
\end{equation}
Now, for any reals $(c,d)$ we take sequences $(a_{n})_{n\geq1}$ and
$(b_{n})_{n\geq1}$ of rational numbers such that $a_{n}\rightarrow c$ and
$b_{n}\rightarrow d$. By the Lebesgue dominated convergence theorem, we pass to
the limit in~(\ref{ab}) (written for $a_{n}$ and $b_{n})$, and obtain that the
equality in~(\ref{ab}) holds for all real numbers. Since the Fourier transform
determines the measure we obtain $(L_{1},L_{2})$ is distributed as
$(N_{1},N_{2})$.
\end{pf}

The next lemma is a step in the proof of Theorem~\ref{char}.

\begin{lemma}
\label{reg}
Let $(\Omega,\mathcal{F},\mathbb{P})$ be a probability space with
$\mathcal{F}$ countably generated, $\mathcal{G}\subset\mathcal{F}$ a sigma
algebra, $Y$ a $\mathcal{G}$-measurable integrable random variable, $X$
integrable and let $g\dvtx \mathbb{R}^{2}\rightarrow\mathbb{R}$ be a continuous and
bounded function. Let $\mathbb{P}^{\omega}$ be a regular version of
$\mathbb{P}$ given $\mathcal{F}$. Then there exists $\Omega_{1}\subset\Omega$
with $\mathbb{P}(\Omega_{1})=1$ such that, for all $\omega\in\Omega_{1}$
\begin{equation}\label{conexpr}
\mathbb{E}^{\omega} \bigl[g \bigl(X,Y(\omega) \bigr) \bigr]=
\mathbb{E}^{\omega} \bigl[g(X,Y) \bigr].
\end{equation}
\end{lemma}

\begin{pf}
It is easy to see that for a simple function $V$ we can find
$\Omega_{V}\subset\Omega$ with $\mathbb{P}(\Omega_{V})=1$ such that, for all
$\omega\in\Omega_{V}$
\[
\mathbb{E}^{\omega} \bigl[g \bigl(X,V(\omega) \bigr) \bigr]=
\mathbb{E}^{\omega} \bigl[g(X,V) \bigr].
\]
Indeed, if $V=\sum_{j=1}^{m}a_{j}I(B_{j})$ with $B_{j}\in\mathcal{G}$ we have
for every $B\in\mathcal{G}$
\begin{eqnarray*}
\mathbb{E} \bigl(I(B)\bigl(g(X,V)|\mathcal{G} \bigr) \bigr)&=& \sum
_{j=1}^{m}\mathbb{E} \bigl(I(B\cap
B_{j})g(X,V) \bigr)
= \sum_{j=1}^{m}\mathbb{E} \bigl(I(B\cap
B_{j})g(X,a_{j}) \bigr)\\
&=& \mathbb{E} \Biggl(I(B)\sum
_{j=1}^{m}I(B_{j})\mathbb{E}
\bigl(g(X,a_{j})|\mathcal{G \bigr) \Biggr)}.
\end{eqnarray*}

Let $V_{n}\rightarrow Y$ a sequence of simple functions. Then we can find a
set $\Omega_{1}\subset\Omega$ with $\mathbb{P}(\Omega_{1})=1$, namely
$\Omega_{1}=\bigcap_{n}\Omega_{V_{n}}$, such that for all $\omega\in\Omega_{1}$
\[
\mathbb{E}^{\omega} \bigl[g \bigl(X,V_{n}(\omega) \bigr) \bigr]=
\mathbb{E}^{\omega} \bigl[g(X,V_{n}) \bigr].
\]
Now, for $\omega$ fixed in $\Omega_{1}$, by Lesbegue dominated convergence
theorem we get (\ref{conexpr}) by passing to the limit.
\end{pf}


\section*{Acknowledgements}

The authors would like to thank the reviewers for their valuable comments and
suggestions which improved the presentation of the paper.
Magda Peligrad was  supported in part by a
Charles Phelps Taft Memorial Fund grant and the NSF  Grant DMS-1208237.




\printhistory
\end{document}